\documentclass[fleqn]{mat01}
\usepackage{times,mathtimy,amssymb,latexsym}%,krepsf,rangecite}
\begin{document}

\setcounter{page}{137}
\firstpage{137}

\def\d{\mbox{\rm d}}
\def\e{\mbox{\rm e}}

\font\xxxxx=msam10 at 10pt
\def\blacksquare{\mbox{\xxxxx{\char'245\ }}}

\newtheorem{theore}{Theorem}
\renewcommand\thetheore{\arabic{section}.\arabic{theore}}
\newtheorem{theor}{\bf Theorem}
\newtheorem{lem}[theor]{\it Lemma}
\newtheorem{rem}[theor]{Remark}

\def\claim{\trivlist \item[\hskip \labelsep{\it Claim.}]}

\newtheorem{pot}[theor]{\it Proof of Theorem}
\renewcommand\thepot{\rm\arabic{pot}}

\title{On Nyman, Beurling and Baez-Duarte's Hilbert space reformulation of the
Riemann hypothesis}

\markboth{Bhaskar Bagchi}{Hilbert space reformulation of the
Riemann hypothesis}

\author{BHASKAR BAGCHI}

\address{Indian Statistical Institute, R.V.~College Post,
Bangalore~560~059, India\\
\noindent E-mail: bbagchi@isibang.ac.in}

\volume{116}

\mon{May}

\parts{2}

\pubyear{2006}

\Date{MS received 21 May 2005; revised 29 April 2006}

\begin{abstract}
There has been a surge of interest of late in an old result of
Nyman and Beurling giving a Hilbert space formulation of the
Riemann hypothesis. Many authors have contributed to this circle
of ideas, culminating in a beautiful refinement due to
Baez-Duarte. The purpose of this little survey is to dis-entangle
the resulting web of complications, and reveal the essential
simplicity of the main results.
\end{abstract}


\maketitle

%\section{Introduction}

\noindent Let$\mathcal{\ H}$\ denote the weighted$\ \,l^{2}$-space
consisting of all sequences $a = \left\{ a_{n}\hbox{:}\
n\in\mathbb{N}\right\} $ of complex numbers such that
$\sum_{n=1}^{\infty}\frac{|a_{n}|^{2}}{n(n+1)}<\infty$. For any
two vectors $\ a,b\in\mathcal{H}$, their inner product is given by
$\left\langle a,b\right\rangle =\sum_{n=1}^{\infty}\frac
{a_{n}\overline{b_{n}}}{n(n+1)}$.  Notice that all bounded
sequences of complex numbers are vectors in this Hilbert space.
For $l=1,2,3,\dots$ let $\ \gamma_{l}\in\mathcal{H}$ be the
sequence
\begin{equation*}
\gamma_{l}=\left\{  \left\{  \frac{n}{l}\right\}\hbox{:}\ n=1,2,3,\dots\right\}.
\end{equation*}
(Here and in what follows,\ $\left\{  x\right\}$ is the fractional part of
a real number $x$.) Also, let $\gamma\in\mathcal{H}$ denote the constant
sequence \
\begin{equation*}
\gamma=\{1,1,1,\dots\}.
\end{equation*}
Recall that a set $A$ of vectors in a Hilbert space $\mathcal{H}$
is said to be {\it total} if the set of all finite linear
combinations of elements of\ $A$\ is dense in $\mathcal{\ H}$,
i.e., if no proper closed subspace of the Hilbert space contains
the set $A$. In terms of these few notions and notations, the
recent result of Baez-Duarte from \cite{Duarte}\ can be given the
following dramatic formulation.

\begin{theor}[\!]
\label{dramatic} The following statements are equivalent{\rm :}

\begin{enumerate}
\renewcommand\labelenumi{\rm (\roman{enumi})}
\leftskip .3pc
\item The Riemann hypothesis{\rm ,}

\item $\gamma$ belongs to the closed linear span of $\ \{\gamma
_{l}\hbox{\rm :}\ l=1,2,3,\dots\}${\rm ,} and

\item the set $\{\gamma_{l}\hbox{\rm :}\ l=1,2,3,\dots\}$ is total in
$\mathcal{\ H}$.\vspace{-.7pc}
\end{enumerate}
\end{theor}

We hasten to add that this is not the statement that the reader
will see in Baez-Duarte's paper. For one thing, the
implications~(ii)~$\Longrightarrow$~(iii) and
(iii)~$\Longrightarrow$~(i) are not mentioned in this paper:
perhaps the author thinks of them as `well-known to experts'. (In
such contexts, an expert is usually defined to be a person who has
the relevant piece of information.) Moreover, the main result in
\cite{Duarte} is not the implication~(i)~$\Longrightarrow$~(ii)
itself, but a `unitarily equivalent' version thereof. More
precisely, the result actually proved in \cite{Duarte} is the
implication~(i)~$\Longrightarrow$~(ii) of Theorem~\ref{refined}
below. In fact, we could not locate in the existing literature the
statement~(iii) of Theorem~\ref{dramatic} (equivalently, of
Theorem~\ref{refined}) as a reformulation of the \hbox{Riemann}
hypothesis. This result may be new. It reveals the Riemann
hypothesis as a version of the central theme of harmonic analysis:
that more or less arbitrary sequences (subject to mild growth
restrictions) can be arbitrarily well approximated by
superpositions of a class of simple periodic sequences (in this
instance, the sequences $\gamma_{l}$).

A second point worth noting is that the particular weight sequence $\big\{\frac
{1}{n(n+1)}\big\}$ used above is not crucial for the validity of Theorem~\ref{dramatic} (though this is the sequence which occurs naturally in its
proof). Indeed, any weight sequence $\{w_{n}\hbox{:}\ n=1,2,3,\dots\}$ satisfying
$\frac{c_{1}}{n^{2}}\leq w_{n}\leq\frac{c_{2}}{n^{2}}$ for all $n$ (for
constants $0<c_{1}\leq c_{2}$) would serve equally well. This is because the
identity map is an invertible linear operator (hence carrying total sets to
total sets) between any two of these weighted $l^{2}$-spaces.

In what follows, we shall adopt the standard practice (in analytic
number theory) of denoting a complex variable by $\ s=\sigma+it$.
Thus $\sigma$ and $t$ are the real and imaginary parts of the
complex number $s$. Recall that {\it Riemann's zeta function} is
the analytic function defined on the half-plane $\{\sigma>1\}$ by
the absolutely convergent series $\zeta
(s)=\sum_{n=1}^{\infty}n^{-s}$. The completed zeta function \
$\zeta^{\ast}$ is \ defined on this half-plane by
$\zeta^{\ast}(s)=\pi^{-s/2}\Gamma (s/2)\zeta(s)$, where $\Gamma$
is Euler's gamma function. As Riemann discovered, $\zeta^{\ast}$
has a meromorphic continuation to the entire complex plane with
only two (simple) poles: at $s=0$ and at $s=1$. Further, it
satisfies the functional equation
$\zeta^{\ast}(1-s)=\zeta^{\ast}(s)$ for all $\ s$. Since $\Gamma$
has poles at the non-positive integers (and nowhere else), it
follows that $\zeta$ has trivial zeros at the negative even
integers. Further, since $\zeta$ is real-valued on the real line,
its zeros occur in conjugate pairs. This trivial observation,
along with the (highly non-trivial) functional equation, shows
that the non-trivial zeros of the zeta function are symmetrically
situated about\break the so-called {\it critical line} $\big
\{\sigma=\frac{1}{2} \big\}$. {\it The Riemann hypothesis} (RH)
conjectures that all these non-trivial zeros actually lie on the
critical line. In view of the symmetry mentioned above, this
amounts to the conjecture that $\zeta$ has no zeros on the
half-plane
\begin{equation*}
\Omega=\left\{s=\sigma+it\hbox{:}\ \sigma>\frac{1}{2}, \ \ -
\infty<t<\infty \right\}.
\end{equation*}
In other words, the Riemann hypothesis is the statement that $\frac{1}{\zeta}$
is analytic on the half-plane $\Omega$. This is the formulation of RH that we
use in this article. Throughout this article, $\Omega$ stands for the
half-plane $\big \{\sigma>\frac{1}{2} \big\}$.

Baez-Duarte's theorem refines an earlier result of the same type
(Theorem~\ref{original} below) proved by Nyman and Beurling
(cf.~\cite{Nyman} and \cite{Beurling}).  Our intention in this
article is to point out that the entire gamut of these results is
best seen inside the {\it Hardy space} $H^{2}(\Omega)$. Recall
that this is the Hilbert space of all analytic functions $F$ on
$\Omega$ such that
\begin{equation*}
\|F\|^{2}\,\hbox{:}\!=\sup_{\sigma>\frac{1}{2}}\frac{1}{2\pi}
\int_{-\infty}^{\infty}\left|F(\sigma+it)\right|^{2} \d t<\infty .
\end{equation*}
It is known that any $F\in H^{2}(\Omega)$ has, almost everywhere on the
critical line, a non-tangential boundary value $F^{\ast}$ such that
\begin{equation*}
\| F \|  ^{2}=\frac{1}{2\pi}\int_{-\infty}^{\infty}\left|
F^{\ast}\left(\frac{1}{2}+it\right)\right|^{2} \d t.
\end{equation*}
Thus $H^{2}(\Omega)$ may be identified (via the isometric
embedding $F\mapsto F^{\ast}$) with a closed subspace of the
$L^{2}$-space of the critical line with respect to the Lebesgue
measure scaled by the factor $\frac{1}{2\pi}$. (This scaling is to
ensure that the Mellin transform $\digamma$, defined while proving
Theorem~\ref{original} below, is an isometry.)

For $0\leq\lambda\leq1$, let $\ F_{\lambda}\in H^{2}(\Omega)$ be defined by
\begin{equation*}
F_{\lambda}(s)=(\lambda^{s}-\lambda)\frac{\zeta(s)}{s},\quad s\in\Omega.
\end{equation*}
Notice that the zero of the first factor at \ $s=1$ cancels the pole of the
second factor, so that $F_{\lambda}$, thus defined, is analytic on $\Omega$.
Also, in view of the well-known elementary estimate (see~\cite{Titchmarsh})
\begin{equation*}
\zeta(s)=O(\left|  s\right|^{\frac{1}{6}}\log\left|  s\right|
),\quad s\in\bar{\Omega}, \quad s\longrightarrow\infty,
\end{equation*}
the factor $\frac{1}{s}$ ensures that $F_{\lambda}\in H^{2}(\Omega)$ for
$0\leq\lambda\leq1$. (Note that, in order to arrive at this conclusion, any
exponent $<\frac{1}{2}$ in the above zeta estimate would have sufficed. But
the exponent $\frac{1}{6}$ happens to be the simplest non-trivial estimate
which occurs in the theory of the Riemann zeta function.) Indeed, under
Riemann hypothesis we have the stronger estimate (Lindelof hypothesis)
\begin{equation}
\zeta(s)=O(|s|^{\epsilon})\quad \hbox{as}\, |s|\longrightarrow\infty,
\hbox{uniformly for}\,s\in\bar{\Omega}, \label{Lindelof}
\end{equation}
for each $\epsilon>0$. (More precisely, under RH, this estimate
holds uniformly on the complement of any given neighbourhood of
$1$ in $\bar{\Omega}$.)

Finally, for $l=1,2,3,\dots$, let $\ G_{l}\in H^{2}(\Omega)$ be defined by
$G_{l}=F_{\frac{1}{l}}$. Thus,
\begin{equation*}
G_{l}(s)=(l^{-s}-l^{-1})\frac{\zeta(s)}{s},\quad s\in\Omega.
\end{equation*}
Also, let $\ E\in H^{2}(\Omega)$ be defined by
\begin{equation*}
E(s)=\frac{1}{s},\quad s\in\Omega.
\end{equation*}
In terms of these notations, the most naural formulation of the
Nyman--Beurling--Baez-Duarte theorem is the following:

\begin{theor}[\!]\label{natural}
The following statements are equivalent{\rm :}

\begin{enumerate}
\renewcommand\labelenumi{\rm (\roman{enumi})}
\leftskip .3pc
\item The Riemann hypothesis{\rm ,}

\item $E$ belongs to the closed linear span of the set $\{G_{l}\hbox{\rm :}\ l=1,2,3,\dots\}${\rm ,} and

\item $E$ belongs to the closed linear span of the set
$\{F_{\lambda}\hbox{\rm :}\ 0\leq\lambda\leq1\}$.
\end{enumerate}
\end{theor}

The plan of the proof is to
verify~(i)~$\Longrightarrow$~(ii)~$\Longrightarrow$~(iii)~$\Longrightarrow$~(i). As we shall see in a little while, except for the
first implication~((i)~$\Longrightarrow$~(ii)), all these implications are
fairly straight forward. In order to prove~(i)~$\Longrightarrow$~(ii), we need
recall that on the half-plane $\{\sigma>1\}, \frac{1}{\zeta}$ is represented
by an absolutely convergent Dirichlet series
\begin{equation}
\sum\limits_{l=1}^{\infty}\mu(l)l^{-s}=\frac{1}{\zeta(s)}.
\end{equation}
Here $\mu(\cdot)$ is the Mobius function. (To determine its
formula, we may formally multiply this Dirichlet series by that of
$\zeta(s)$ and equate coefficients to get the recurrence relation
$\sum_{l|n} \mu(l)=\delta_{1n}$. Solving this, one can show that
$\mu(\cdot)$ takes values in $\{0,+1,-1 \}$ and hence the
Dirichlet series for $\frac {1}{\zeta}$ is absolutely convergent
on $\{\sigma > 1 \}$. Indeed, $\mu(l) = 0$ if $l$ has a repeated
prime factor, $\mu(l) = +1$ if $l$ has an even number of distinct
prime factors, and $\mu(l) = -1$ if $l$ has an odd number of
distinct prime factors. But, for our limited purposes, all this is
unnecessary.) What we need is an old theorem of Littlewood
(see~\cite{Titchmarsh}) to the effect that for the validity of the
Riemann hypothesis, it is necessary (and sufficient) that the
Dirichlet series displayed above converges uniformly on compact
subsets of $\Omega$. Actually, we need the following quantitative
version of this theorem of Littlewood.

%We shall use the following more or less well-known lemma.%%

\begin{lem}\label{Littlewood}
If the Riemann hypothesis holds then for each $\epsilon>0$ and
each $\delta>0${\rm ,} we have
$\sum_{l=1}^{L}\mu(l)l^{-s}=O((\left| t\right|  +1)^{\delta})$
uniformly for $L=1,2,3,\dots$ and uniformly for $s=\sigma+it$ in
the half-plane $\big\{ \sigma>\frac{1}{2}+\epsilon \big\}$. {\rm
(}Thus the implied constant depends only on $\epsilon$ and
$\delta$.{\rm )}
\end{lem}

Since
Lemma~\ref{Littlewood} is more or less well-known, we omit its
proof. It may be proved by a minor variation in the original
proof of Littlewood's theorem quoted above. (Note that, with the
aid of a little `normal family' argument, Littlewood's theorem
itself is an easy consequence of this lemma.)

\setcounter{theor}{1}
\begin{pot}{\rm
(i)~$\Rightarrow$~(ii). Assume RH. For
positive integers $L$ \ and any small real number $\epsilon>0$, let
$H_{L,\epsilon}\in H^{2}(\Omega)$ be defined by
\begin{equation*}
H_{L,\epsilon}=\sum\limits_{l=1}^{L}\frac{\mu(l)}{l^{\epsilon}}G_{l}.
\end{equation*}
Thus each $H_{L,\epsilon}$ is in the linear span of $\{G_{l}\hbox{:}\ l\geq1\}$. Note that
\begin{equation*}
H_{L,\epsilon}(s)=\frac{\zeta(s)}{s}\left(\sum\limits_{l=1}^{L}\frac{\mu
(l)}{l^{s+\epsilon}}-\sum\limits_{l=1}^{L}\frac{\mu(l)}{l^{1+\epsilon}}\right),\quad
s\in\bar{\Omega}.
\end{equation*}
Therefore, by the theorem of Littlewood quoted above, for any fixed
$\epsilon>0$,
\begin{equation*}
H_{L,\epsilon}(s)\longrightarrow H_{\epsilon}(s)\quad \hbox{for}\,s \ \hbox{in
the critical line, as}\,L\longrightarrow\infty.
\end{equation*}
Here,
\begin{equation*}
H_{\epsilon}(s)\,\hbox{:}\!=\frac{\zeta(s)}{s}\left(\frac{1}{\zeta(s+\epsilon)}-\frac{1}
{\zeta(1+\epsilon)}\right).
\end{equation*}
Also, by the estimate~(\ref{Lindelof}) and Lemma~\ref{Littlewood},
$H_{L,\epsilon}$ is bounded by an absolutely square integrable
function on the critical line. Therefore, by Lebesgue's dominated
convergence theorem, we have, for each fixed $\epsilon>0$,
\begin{equation*}
H_{L,\epsilon}\longrightarrow H_{\epsilon} \quad \hbox{in the norm of}
\ H^{2}(\Omega) \ \hbox{as} \ L\longrightarrow\infty.
\end{equation*}
Since $H_{L,\epsilon}$ is in the linear span of $\{G_{l}\hbox{:}\ l=1,2,3,\dots\},$ it
follows that, for each $\epsilon>0$, $H_{\epsilon}$ is in the closed linear
span of $\{G_{l}\hbox{:}\ l=1,2,3,\dots\}$. Now note that, since $\zeta$ has a pole at
$s=1$,
\begin{equation*}
H_{\epsilon}(s)\longrightarrow\frac{1}{s}=E(s)\ \hbox{for} \ s\
\hbox{in the critical line, as} \ \epsilon\searrow 0.
\end{equation*}
Therefore, in order to show that $E$ is in the closed linear span of
$\{G_{l}\hbox{:}\ l=1,2,3,\dots\}$ and thus complete this part of the proof, it
suffices to show that $H_{\epsilon},\,0<\epsilon<\frac{1}{2},$ are uniformly
bounded in modulus on the critical line by an absolutely square integrable
function. Then, another application of Lebesgue's dominated convergence theorem would yield
\begin{equation*}
H_{\epsilon}\longrightarrow E \ \hbox{in the norm of} \ H^{2}(\Omega) \ \hbox{as} \ \epsilon\searrow 0.
\end{equation*}
Consider the entire function $\xi(s)\,\hbox{:}\!=s(1-s)\zeta^{\ast}(s)=s(1-s)\pi
^{-s/2}\Gamma\big(\frac{s}{2}\big)\zeta(s)$. It has the Hadamard factorisation
\begin{equation*}
\xi(s)=\xi(0)\prod\limits_{\rho}\left(1-\frac{s}{\rho}\right),
\end{equation*}
where the product is over all the non-trivial zeros $\rho$ of the
Riemann zeta function. This product converges provided the zeros
$\rho$ and $1-\rho$ are grouped together. In consequence, with a
similar bracketing, we have
\begin{equation*}
|\xi(s)|=|\xi(0)|\prod\limits_{\rho}\left| 1-\frac{s}{\rho}\right|.
\end{equation*}
Now, under RH, each $\rho$ has real part $=\frac{1}{2}$.
Therefore, for $s$ in the closed half-plane $\bar{\Omega}$, we
have $\big|1-\frac{s}{\rho}
\big|\leq\big|1-\frac{s+\epsilon}{\rho}\big|$. Multiplying this
trivial inequality over all $\rho$, we get
\begin{equation*}
|\xi(s)|\leq|\xi(s+\epsilon)|, \ \ s\in\bar{\Omega},\,\epsilon>0.
\end{equation*}
(Aside: conversely, the above inequality clearly implies RH. Thus,
this simple looking inequality is a reformulation of RH.) In other
words, we have, for $s\in\bar{\Omega}$,
\begin{align*}
\left|  \frac{\zeta(s)}{\zeta(s+\epsilon)}\right|\! \leq\pi^{-\epsilon
/2}\left|  \frac{(s+\epsilon)(1-\epsilon-s)}{s(1-s)}\right|  \left|
\frac{\Gamma((s+\epsilon)/2)}{\Gamma(s/2)}\right|\! \leq c\left|  \frac
{\Gamma((s+\epsilon)/2)}{\Gamma(s/2)}\right|
\end{align*}
for some absolute constant $c>0$. But, by Sterling's formula (see
\cite{Lang} for instance), the gamma ratio on the extreme right is
bounded by constant times $|s|^{\epsilon/2}$, uniformly for
$s\in\bar{\Omega}$. Therefore we get
\begin{equation*}
\left|  \frac{\zeta(s)}{\zeta(s+\epsilon)}\right|  \leq
c|s|^{\epsilon /2},\quad s\in\bar{\Omega},
\end{equation*}
for some other absolute constant $c>0$. In conjunction with the
estimate~(\ref{Lindelof}), this implies
\begin{equation*}
|H_{\epsilon}(s)|\leq c|s|^{-3/4},\quad s\in\bar{\Omega},
\end{equation*}
for $0<\epsilon<\frac{1}{2}$. Since $s\longmapsto$ $c|s|^{-3/4}$ is square
integrable on the critical line, we are done. This proves the
implication~(i)~$\Rightarrow$~(ii).

Since $\{G_{l}\hbox{:}\ l=1,2,3,\dots\}\subseteq\{F_{\lambda}\hbox{:}\ 0\leq\lambda\leq1\}$, the
implication~(ii)~$\Rightarrow$~(iii) is trivial. To prove~(iii)~$\Rightarrow$~(i), suppose RH is false. Then there is a zeta-zero $\rho\in\Omega$. Since
$\zeta(\rho)=0$, it follows that $F_{\lambda}(\rho)=0$ for all $\lambda
\in(0,1]$. Thus the set $\{F_{\lambda}\hbox{:}\ \lambda\in(0,1]\}$ (and hence also its
closed linear span) is contained in the proper closed subspace $\{F\in
$\ $H^{2}(\Omega)\hbox{:}\ F(\rho)=0\}$ of \ $H^{2}(\Omega)$. (It is a closed subspace
since evaluation at any fixed \ $\rho\in\Omega$ is a continuous linear
functional: $H^{2}(\Omega)$ is a functional Hilbert space.) Since \ $E$
belongs to the closed linear span of this set, it follows that
$0=E(\rho)=\frac{1}{\rho}$. Hence $0=1$: the ultimate contradiction!
This proves~(iii)~$\Longrightarrow$~(i).} \hfill $\blacksquare$
\end{pot}

\setcounter{theor}{3}
\begin{rem}{\rm
\label{First} Since $\mu(l)=0$ unless $l$ is square-free, the
functions $H_{L,\epsilon}$ introduced in the course of the above
proof are in the linear span of the set $\{G_{l}\hbox{:}\ l$
square-free$\}$. Thus, the proof actually shows that RH implies
(and hence is equivalent to) that $E$ belongs to\break the closed
linear span of the thinner set $\{G_{l}\hbox{:}\ l$
square-free$\}$ in $H^{2}(\Omega)$.}
\end{rem}

Now let $L^{2}((0,1])$ be the Hilbert space of complex-valued absolutely
square integrable functions (modulo almost everywhere equality) on the
interval $(0,1]$. For $0\leq\lambda\leq1$, let $f_{\lambda}\in L^{2}((0,1])$
be defined by
\begin{equation*}
f_{\lambda}(x)=\left\{\frac{\lambda}{x}\right\}-
\lambda\left\{\frac{1}{x}\right\}, \quad x\in(0,1].
\end{equation*}
(Recall that $\{ \cdot\}$ stands for the fractional part.)\ Let $\mathbf{1}\in$
$L^{2}((0,1])$ denote the constant function $=1$ on $(0,1]$. Thus,
\begin{equation*}
\mathbf{1(}x)=1\mathbf{,} \ x\in(0,1].
\end{equation*}
In terms of these notations, the original theorem of Nyman and
Beurling may be stated as follows:

\begin{theor}[\!]\label{original}
The following statements are equivalent\hbox{\rm :}

\begin{enumerate}
\renewcommand\labelenumi{\rm (\roman{enumi})}
\leftskip .35pc
\item The Riemann hypothesis{\rm ,}

\item $\mathbf{1}$ is in the closed linear span in $L^{2}((0,1])$
of the set $\{f_{\lambda}\hbox{\rm :}\ 0\leq\lambda\leq1\}${\rm ,}

\item the set $\{f_{\lambda}\hbox{\rm :}\ 0\leq\lambda\leq1\}$ is total in
$L^{2}((0,1])$.
\end{enumerate}
\end{theor}

\begin{proof}
One defines the \emph{Fourier--Mellin transform \
}$\digamma\hbox{:}\ L^{2} ((0,1])\longrightarrow H^{2}(\Omega)$ by
\begin{equation}
\digamma(f)(s)=\int_{0}^{\infty}x^{s-1}f(x)\d x,\quad
s\in\Omega,  \ f\in L^{2}((0,1]).
\end{equation}
It is well-known that $\digamma$, thus defined, is an isometry.
For completeness, we sketch a proof. Since
$s\longmapsto(x\longmapsto x^{s-1})$ is an $L^{2}((0,1])$-valued
analytic function on $\Omega$, it follows that $\digamma(f)$ is
analytic on $\Omega$ for each $f\in L^{2}((0,1])$. For
$\lambda\in[0,1]$, let $\Psi_{\lambda}\in L^{2}((0,1])$ denote the
indicator function of the interval $(0,\lambda)$. Using the
well-known identity
\begin{equation*}
\frac{1}{\pi}\int_{-\infty}^{+\infty}\frac{\e^{iux}}{1+x^{2}}
\d x=\e^{-|u|},\quad u\in\mathbb{R},
\end{equation*}
one sees that $\|\digamma(\Psi_{\lambda})\| ^{2}=\|
\Psi_{\lambda}\|  ^{2}<\infty$~--~hence
$\digamma(\Psi_{\lambda})\in H^{2}(\Omega)$~--~and, more
generally, $\|\digamma(\Psi_{\lambda })-\digamma(\Psi_{\mu})\|
^{2}=\| \Psi_{\lambda}-\Psi_{\mu }\|  ^{2}$ for
$\lambda,\mu\in[0,1]$. Since $\{\Psi_{\lambda}\hbox{:}\
\lambda\in[0,1]\}$ is a total subset of $L^{2}((0,1])$, this
implies that $\digamma$ maps $L^{2}((0,1])$ isometrically into
$H^{2}(\Omega)$.

We begin with a computation of the Melin transform of
$f_{\lambda}$.

\begin{claim}
\begin{equation}
\digamma(f_{\lambda})=-F_{_{\lambda}},\quad 0\leq\lambda\leq1.
\end{equation}
To verify this claim, begin with $s=\sigma+it$, \ $\sigma>1$.
Then,
\begin{align*}
\int_{0}^{1}\left\lbrace\frac{\lambda}{x}\right\rbrace
x^{s-1}\d x&=\lambda\int_{0}^{1}x^{s-2}\d
x-\int_{0}^{1}\left\lfloor \frac{\lambda}{x}\right\rfloor
x^{s-1}\d x\\[.5pc]
&=\frac{\lambda}{s-1}-\int_{0}^{1}\left\lfloor \frac{\lambda
}{x}\right\rfloor x^{s-1}\d x.
\end{align*}
But,
\begin{align*}
\int_{0}^{1}\left\lfloor \frac{\lambda}{x}\right\rfloor
x^{s-1}\d x &=\sum_{n=1}^{\infty}n
\int_{\lambda/(n+1)}^{\lambda/n}x^{s-1}\d x\\[.5pc]
&=\frac{\lambda^{s}}{s}\sum\limits_{n=1}^{\infty}n(n^{-s}-(n+1)^{-s}).
\end{align*}
Now, the partial sum $\sum_{n=1}^{N}n(n^{-s}-(n+1)^{-s})$
telescopes to $-N(N+1)^{-s}+$ $\sum_{n=1}^{N}n^{-s}$. Since
$\sigma>1$, letting $N\longrightarrow\infty$, we get
$\sum_{n=1}^{\infty}n(n^{-s} -(n+1)^{-s})=\zeta(s)$. Thus,
\begin{equation*}
\int_{0}^{1}\left\lbrace\frac{\lambda}{x}\right\rbrace
x^{s-1}\d x=\frac{\lambda}{s-1} -\lambda^{s}\frac{\zeta(s)}{s}.
\end{equation*}
In particular, taking $\lambda=1$ here, one gets
\begin{equation*}
\int_{0}^{1}\left\lbrace\frac{1}{x}\right\rbrace x^{s-1}\d
x=\frac{1}{s-1}-\frac{\zeta(s)}{s}.
\end{equation*}
Multiplying the second equation by $\lambda$ and
subtracting the result from the first, we arrive at
\begin{equation*}
\int_{0}^{1}f_{\lambda}(x)x^{s-1}\d x
=-(\lambda^{s}-\lambda)\frac {\zeta(s)}{s}=-F_{\lambda}(s)
\end{equation*}
for $s$ in the half-plane $\{\sigma>1\}$. Since both sides of this
equation are analytic in the bigger half-plane $\Omega$, this
equation continues to hold for $s\in\Omega$. This proves the
Claim.\vspace{.5pc}
\end{claim}

\noindent (i)~$\Longrightarrow$~(ii). Assume RH. Then, by
Theorem~\ref{natural}, $E=\digamma(\mathbf{1})$ belongs to the
closed linear span of $\{F_{\lambda
}=-\digamma(f_{\lambda})\hbox{:}\ 0\leq\lambda\leq1\}$. Since
$\digamma$ is an isometry, this shows that $\mathbf{1}$ belongs to
the closed linear span of the set $\{f_{\lambda}\hbox{:}\
0\leq\lambda\leq1\}$. Thus (i)~$\Longrightarrow$~(ii).\vspace{.5pc}

\noindent (ii)~$\Longrightarrow$~(iii). Let $\mathbf{1}$ be in the closed
linear span in $L^{2}((0,1])$ of $\{f_{\lambda}\hbox{:}\
0\leq\lambda\leq1\}$. Applying $\digamma$, it follows that $E$ is
in the closed linear span (say $\mathcal{N}$) of
$\{F_{\lambda}\hbox{:}\ 0\leq\lambda\leq1\}$. For $\mu\in(0,1]$,
let $\Theta_{\mu}\in H^{\infty}(\Omega)$ (the Banach algebra of
bounded analytic functions on $\Omega$) be defined by
\begin{equation*}
\Theta_{\mu}(s)=\mu^{s-\frac{1}{2}},\quad s\in\Omega.
\end{equation*}
We have $|\Theta_{\mu}(s)|=1$ for $s$ in the critical line. That
is, $\Theta_{\mu}$ is an inner function. In consequence, the
linear operators $M_{\mu}\hbox{:}\ H^{2}(\Omega)\longrightarrow
H^{2}(\Omega)$ defined by
\begin{equation*}
M_{\mu}(F)=\Theta_{\mu}F \ \hbox{(point-wise product),} \ F\in
H^{2}(\Omega),
\end{equation*}
are isometries. (Since
$\Theta_{\lambda}\Theta_{\mu}=\Theta_{\lambda\mu}$, it follows
that $M_{\lambda}M_{\mu}=M_{\lambda\mu}$ for
$\lambda,\mu\in(0,1]$. Thus $\{M_{\mu}\hbox{:}\ \mu\in(0,1]\}$ is
a semi-group of isometries on $H^{2} (\Omega)$ modelled after the
multiplicative semi-group $(0,1]$.) Trivially, for
$0\leq\lambda\leq1$ and $0<\mu\leq1$, we have
\begin{equation*}
M_{\mu}(F_{\lambda})=\Theta_{\mu}F_{\lambda}=\mu^{-1/2}
(F_{\lambda\mu}-\lambda F_{\mu}).
\end{equation*}
This shows that the closed subspace $\mathcal{N}$ spanned by the
$F_{\lambda}$'s is invariant under the semi-group
$\{M_{\mu}\hbox{:}\ \mu\in(0,1]\}$:
\begin{equation*}
M_{\mu}(\mathcal{N)}\subseteq\mathcal{N},\quad\mu\in(0,1].
\end{equation*}
Since $E\in\mathcal{N}$, it follows that
$M_{\mu}(E)\in\mathcal{N}$ for $\mu\in(0,1]$. But we have the
trivial computation
\begin{equation*}
\digamma(\Psi_{\lambda})=\lambda^{1/2}M_{\lambda}(E),\quad
0<\lambda\leq 1.
\end{equation*}
Thus, $\{\digamma(\Psi_{\lambda})\hbox{:}\ 0\leq\lambda\leq1\}$ is
contained in the closed linear span $\mathcal{N}$ of
$\{\digamma(f_{\lambda})\hbox{:}\ 0\leq \lambda\leq1\}$. Since
$\digamma$ is an isometry, it follows that
$\{\Psi_{\lambda}\hbox{:}\ 0\leq\lambda\leq1\}$ is contained in
the closed linear span in $L^{2}((0,1])$ of the set
$\{f_{\lambda}\hbox{:}\ 0\leq\lambda\leq1\}$. Since the first set
is clearly total in $L^{2}((0,1])$, it follows that so is the
second. Thus (ii)~$\Longrightarrow$~(iii).\vspace{.4pc}

\noindent (iii)~$\Longrightarrow$~(i). Clearly (iii) implies that the closed
linear span of $\{f_{\lambda}\hbox{:}\ 0\leq\lambda\leq1\}$
contains $\mathbf{1}$ and hence, applying $\digamma$, the closed
linear span of $\{F_{\lambda}\hbox{:}\  0\leq \lambda\leq1\}$
contains $E$. Therefore, by Theorem~\ref{natural}, Riemann
hypothesis follows. Thus (iii)~$\Longrightarrow$~(i). \hfill
$\blacksquare$
\end{proof}

\begin{rem}\label{second}
{\rm It is instructive to compare the proof of
Theorem~\ref{original} with Beurling's original proof as given in
\cite{Donoghue}. Our proof makes it clear that the heart of the
matter is very simple: Riemann hypothesis amounts to the existence
of approximate inverses to the zeta function in a suitable
function space (viz. the weighted Hardy space of analytic
functions on $\Omega$ with the weight function $|E(s)|^{2}$). The
simplification in its proof is achieved by \hbox{Baez-Duarte's}
perfectly natural and yet vastly illuminating observation that,
under RH, these approximate inverses are provided by the partial
sums of the Dirichlet series for $\frac{1}{\zeta}$. In contrast,
Beurling's original proof is a clever and ill-motivated
application of Phragmen--Lindelof type arguments. (We have not
seen Nyman's original proof.) To be fair, we should however point
out that such arguments are now hidden under the carpet: they
occur in the proofs (not presented here) of the conditional
estimate~(\ref{Lindelof}) and Lemma~\ref{Littlewood}.}
\end{rem}

Let $\mathcal{M}$ be the closed subspace of
$L^{2}((0,1])$ consisting of the functions which are almost
everywhere constant on each of the sub-intervals
$\big(\frac{1}{n+1},\frac{1}{n}\big],\;n=1,2,3,\dots$. Since each
element of $\mathcal{M}$ is almost everywhere equal to a unique
function which is everywhere constant on these sub-intervals, we
may (and do) think of $\mathcal{M}$ as the space of all such
(genuine) piece-wise constant functions. As a closed subspace of a
Hilbert space, $\mathcal{M}$ is a Hilbert space in its own right.

For $l=1,2,3,\dots$, let $g_{l}\in$\ $L^{2}((0,1])$ be defined by
\begin{equation*}
g_{l}(x)=\left\lbrace\frac{1}{lx}\right\rbrace
-\frac{1}{l}\left\lbrace\frac{1}{x}\right\rbrace, \quad x\in(0,1].
\end{equation*}
Thus, $g_{l}=f_{1/l}$, $l=1,2,3,\dots$.

Notice that we have $g_{l}(x)=\frac{1}{l}\left\lfloor
\frac{1}{x}\right\rfloor -\left\lfloor \frac{1}{lx}\right\rfloor$.
Also, for $x\in\big(\frac{1}{n+1}, \frac{1}{n}\big], n =
1,2,3,\dots$, $\frac{1}{lx}\in$
$\big[\frac{n}{l},\frac{n+1}{l}\big)$, and no integer can be in
the interior of the latter interval, so that $\left\lfloor
\frac{1}{lx}\right\rfloor =\left\lfloor \frac{n}{l}\right\rfloor$;
also, $\left\lfloor \frac{1}{x}\right\rfloor =n$ for
$x\in\big(\frac{1}{n+1},\frac{1}{n}\big]$. Thus we get
\begin{equation}
g_{l}(x)=g_{l}\left(\frac{1}{n}\right)=\left\lbrace\frac{n}{l}\right\rbrace
,\quad x\in\left(\frac{1}{n+1},\frac {1}{n}\right].
\end{equation}
In consequence,
\begin{equation*}
g_{l}\in\mathcal{M},\quad l=1,2,3,\dots.
\end{equation*}
The refinement due to Baez-Duarte of the Beurling--Nyman theorem
may now be stated as follows. (However, as already stated, the
implication (i)~$\Longrightarrow$~(ii) of this theorem is its only
part which explicitly occurs in \cite{Duarte}.)

\begin{theor}[\!]\label{refined}
The following are equivalent\hbox{\rm :}

\begin{enumerate}
\renewcommand\labelenumi{\rm (\roman{enumi})}
\leftskip .3pc
\item The Riemann hypothesis{\rm ,}

\item $\mathbf{1}$ belongs to the closed linear span of
$\{g_{l}\hbox{\rm :}\ l=1,2,3,\dots\}${\rm ,} and

\item $\{g_{l}\hbox{\rm :}\ l=1,2,3,\dots\}$ is a total set in
$\mathcal{M}$.
\end{enumerate}
\end{theor}

\begin{proof}
Putting $\lambda=\frac{1}{l}$ in the formula~(4), we get
\begin{equation*}
\digamma(g_{l})=-G_{l},\quad l=1,2,3,\dots.
\end{equation*}
Since, under RH, $E=\digamma(\mathbf{1)}$ is in the closed linear
span of $\{G_{l}=-\digamma(g_{l})\hbox{:}\ l=1,2,3,\dots\}$ and
$\digamma$ is an isometry, it follows that $\mathbf{1}$ is in the
closed linear span of $\{g_{l}\hbox{:}\ l=1,2,3,\dots\}$. Thus
(i)~$\Longrightarrow$~(ii).

Now, for positive integers $m$, define the linear operators
$T_{m}\hbox{:}\ \mathcal{M} \longrightarrow \mathcal{M}$ by
\begin{equation*}
(T_{m}f)(x)=\begin{cases}
m^{1/2}f(mx), &\hbox{if} \ x\in\big(0,\frac{1}{m}\big],\\[.5pc]
0,            &\hbox{if} \ x\in\big(\frac{1}{m},1\big].
\end{cases}
\end{equation*}
Clearly each $T_{m}$ is an isometry. (We have
$T_{m}T_{n}=T_{mn}$~--~thus $\{T_{m}\hbox{:}\ m=1,2,3,\dots\}$ is
a semigroup of isometries modelled after the multiplicative
semi-group of positive integers.) Also, it is easy to see that
\begin{equation*}
T_{m}(g_{l})=m^{1/2}\left(g_{lm}-\frac{g_{m}}{l}\right)
\end{equation*}
for any two positive integers $l,m$. Thus the closed linear span
$\mathcal{K}$ of the vectors $g_{l}$, $l=1,2,3,\dots$ is invariant
under this semi-group. Further, letting $\Phi_{n}\in\mathcal{M}$
denote the indicator function of the interval
$\big(0,\frac{1}{n}\big]$, one has
\begin{equation*}
T_{m}(\Phi_{n})=m^{1/2}\Phi_{mn}.
\end{equation*}
Thus, if $\mathcal{K}$ contains $\mathbf{1}=\Phi_{1}$ then it
contains $\Phi_{n}$ for all $n$. Since $\{\Phi_{n}\hbox{:}\
n=1,2,3,\dots\}$ is clearly a total subset of $\mathcal{M}$, it
then follows that $\mathcal{K=M}$, so that $\{g_{l}\hbox{:}\
l=1,2,3,\dots\}$ is a total subset of $\mathcal{M}$. Thus
(ii)~$\Longrightarrow$~(iii).

Lastly, if $\{g_{l}\hbox{:}\ l=1,2,3,\dots\}$ is a total subset of
$\mathcal{M}$ then, in particular its closed linear span contains
$\mathbf{1,}$ and hence the closed linear span of $
\{G_{l}=-\digamma(g_{l})\}$ contains $E=\digamma(\mathbf{1})$, so
that RH follows by Theorem~\ref{natural}. Thus
(iii)~$\Longrightarrow$~(i). \hfill $\blacksquare$
\end{proof}

\setcounter{theor}{0}
\begin{pot}{\rm  Let $U\hbox{:}\
\mathcal{M\longrightarrow H}$ be the unitary defined by
\begin{equation*}
U(f)=\left\lbrace f\left(\frac{1}{n}\right)\hbox{:}\
n=1,2,3,\dots\right\rbrace,\quad f\in\mathcal{M}.
\end{equation*}
Since $U(\mathbf{1})=\gamma$ and (in view of equation~(5))
$U(g_{l})=\gamma_{l}$, this theorem is a straightforward
reformulation of theorem~\ref{refined}.} \hfill $\blacksquare$
\end{pot}

\setcounter{theor}{7}
\begin{rem}
{\rm In view of Remark~\ref{First}, Riemann hypothesis actually
implies (and hence is equivalent to) the statement that $\gamma$
belongs to the closed linear span in $\mathcal{H}$ of the much
thinner set $\{\gamma_{l}\hbox{:}\ l$ square-free$\}$.}
\end{rem}

So where does the undoubtedly elegant reformulation of RH in
Theorem~\ref{dramatic} leave us? One possible approach is as
follows. For positive integers $L$, let $D(L)$ denote the distance
of the vector $\gamma \in\mathcal{H}$ from the $(L-1)$-dimensional
subspace of $\mathcal{H}$ spanned by
$\gamma_{1},\gamma_{2,}\dots,\gamma_{L}$. In view of
Theorem~\ref{dramatic}, RH is equivalent to the statement
$D(L)\longrightarrow0$ as $L\longrightarrow \infty$. So one might
try to estimate $D(L)$. Indeed, as a discrete analogue of a
conjecture of Baez-Duarte {\it et~al} \cite{DBLS}, one might
expect that $D^{2}(L)$ is asymptotically equal to $\frac{A}{\log
L}$ for $A=2+C-\log (4\pi)$, where $C$ is Euler's constant. (But,
of course, this is far stronger than RH itself.) A~standard
formula gives $D^{2}(L)$ as a ratio of two Gram determinants,
i.e., determinants with the inner products $\left\langle
\gamma_{l},\gamma_{m}\right\rangle $ as entries. It is easy to
write down these inner products as finite sums involving the
logarithmic derivative of the gamma function. But such formulae
are hardly suitable for calculation/estimation of determinants. In
any case, it will be a sad day for Mathematics when (and if) the
Riemann hypothesis is proved by a brute-force calculation! Surely
a dramatically new and deep idea is called for. But then, as a
wise man once said, it is fool-hardy to predict~--~specially the
future!


\begin{thebibliography}{9}
\bibitem{Beurling} Beurling~A, A closure problem related to the
Riemann zeta function, {\it Proc. Natl. Acad. Sci.} {\bf 41}
(1955) 312--314

\bibitem{Duarte} Baez-Duarte~L, A strengthening of the Nyman--Beurling
criterion for the Riemann hypothesis, {\it Atti Acad. Naz. Lincei}
{\bf 14} (2003) 5--11

\bibitem{DBLS} Baez-Duarte~L, Balazard~M, Landreau~B and Saias~E,
Notes sur la fonction $\zeta$ de Riemann 3, {\it Advances in
Math.} {\bf 149} (2000) 130--144

\bibitem{Donoghue} Donoghue~Jr~W~F, Distributions and Fourier transforms
(Academic Press) (1969)

\bibitem{Lang} Lang~S, Complex Analysis (Springer Verlag) (1992)

\bibitem{Nyman} Nyman~B, On some groups and semi-groups of translations,
Ph.~D.~Thesis (Uppsala) (1950)

\bibitem{Titchmarsh} Titchmarsh~E~C, The theory of the Riemann zeta
function (Oxford Univ. Press) (1951)
\end{thebibliography}
\end{document}